\documentstyle[twoside,fleqn]{article}
\newcounter{zacountsec}
\newcommand{\eh}{\hfill}\newlength{\sperr}
\newcommand{\Title}[1]{{\large \bf #1}}
\newcommand{\Author}[1]{{\sc #1}}
\newcommand{\Section}[1]{{\stepcounter{zacountsec}\vspace{3mm}%
\hspace*{18mm}\normalsize\bf\arabic{zacountsec}. \parbox[t]{150mm}{ #1 }}}
\newcommand{\Newpara}{\\[0.1cm]\mbox{}\hspace*{9mm}}
\newenvironment{Abstract}{\begin{minipage}[t]{177mm}\em }%
{\end{minipage}}
\newenvironment{thm}[2]{\begin{sloppypar}%
{#1 #2.}\em{}}%
{\end{sloppypar}}
\newcommand{\theo}[3]{\begin{thm}{#1}{#2} #3\end{thm}}
\newcommand{\Theorem}{\hspace*{9mm}{\settowidth{\sperr}{\rm Theorem}%
\parbox[t]{1.3\sperr}{\rm T\eh h\eh e\eh o\eh r\eh e\eh m\eh } }}
\newcommand{\Remark}{\hspace*{9mm}{\settowidth{\sperr}{\rm Remark}%
\parbox[t]{1.3\sperr}{\rm R\eh e\eh m\eh a\eh r\eh k\eh } }}
\newcommand{\Lemma}{\hspace*{9mm}{\settowidth{\sperr}{\rm Lemma}%
\parbox[t]{1.3\sperr}{\rm L\eh e\eh m\eh m\eh a\eh } }}

\newcommand{\proof}{\hspace*{9mm}{\settowidth{\sperr}{\rm Proof}%
\parbox[t]{1.3\sperr}{\rm P\eh r\eh o\eh o\eh f\eh. } }}

\newcommand{\R}{{\rm I}\!{\rm R}}

\newcounter{zalit}

\newenvironment{Acknowledgements}{\vspace{3mm}%
\hspace*{18mm}{\bf Acknowledgements}\\[0.3cm]\begin{minipage}[t]{177mm}%
\small \em}{\end{minipage}}
\newenvironment{References}{%
\Section{References}%
\begin{small}\begin{list}{\arabic{zalit} }{\usecounter{zalit}
\itemsep0mm \parsep0mm\settowidth{\labelwidth}{\small\rm 88}\labelsep0mm
\setlength{\leftmargin}{\labelwidth}}}%
{\end{list}\end{small}}
\newlength{\addro}\newlength{\addrt}
\newenvironment{Address}{{\em Addresses: }%
\begin{minipage}[t]{\addrt}}%
{\end{minipage}}

\newcommand{\bH}{{\bf H}}

\newcommand{\cD}{{\cal D}}

\newcommand{\cH}{{\cal H}}

\newcommand{\cQ}{{\cal Q}}

\newcommand{\an}{\alpha}
\newcommand{\bn}{\beta}

\newcommand{\dn}{\delta}

\newcommand{\ad}{_{\alpha}}

\newcommand{\bd}{_{\beta}}

\newcommand{\aad}{_{\alpha\alpha}}

\newcommand{\abd}{_{\alpha\beta}}
\newcommand{\bad}{_{\beta\alpha}}
\newcommand{\au}{^{(\alpha)}}

\newcommand{\ju}{^{(j)}}
\newcommand{\Sum}{\displaystyle\sum\limits}
\newcommand{\Int}{\int\limits}

\newcommand{\Frac}[2]{\frac{\textstyle #1}{\textstyle #2}}
\newcommand{\Min}[1]{\mathop{{\rm min}}\limits_{#1}}

\newcommand{\Sup}[1]{\mathop{\rm sup}\limits_{#1}}

\newcommand{\reduction}[2]{\left. #1 \right|_{#2}}

\newcommand{\Lim}[1]{\mathop{{\rm lim}}\limits_{#1}}
\newcommand{\diag}{{\rm diag}}
\newcommand{\be}{\begin{equation}}
\newcommand{\ee}{\end{equation}}
\begin{document}
\vspace*{15mm}\hspace*{18mm}
\begin{minipage}[t]{157mm}
\Author{Motovilov, Alexander K.}
\vspace*{0.4cm}

\Title{ Removal of the Resolvent-like
  Dependence on the Spectral Parameter
\newline from Perturbations\footnote{Published 
in Proceedings of the International Congress on
Industrial and Applied Mathematics (Hamburg,
July 1995), Zeitschrift f\"ur Angewandte Mathematik und Mechanik
(ZAMM), Special Issue 2 (ISBN 3-05-501745-5), 229--232 (1996). 
LANL E-print {\tt math.SP/9810044}.} }
\end{minipage}
\vspace*{3.5mm}

\begin{Abstract}
The spectral problem $(A + V(z))\psi=z\psi$ is considered with 
$A$, a self-adjoint operator.  The perturbation $V(z)$  is 
assumed to depend on the spectral parameter $z$ as resolvent of 
another self-adjoint operator $A':$ $V(z)=$ 
\mbox{$-B(A'-z)^{-1}B^{*}$.}  It is supposed that the operator 
$B$ has a finite Hilbert-Schmidt norm  and spectra of the 
operators $A$ and $A'$ are separated.  Conditions are formulated 
when the perturbation $V(z)$ may be replaced with a 
``potential'' $W$ independent of $z$ and such that the operator 
$H=A+W$ has the same spectrum and the same eigenfunctions (more 
precisely, a part of spectrum and a respective part of 
eigenfunctions system) as the initial spectral problem.  The 
operator $H$ is constructed as a solution of the non--linear 
operator equation $H=A+V(H)$ with a specially chosen 
operator--valued function $V(H)$. In the case if the initial 
spectral problem corresponds to a two--channel variant of the 
Friedrichs model, a basis property of the eigenfunction system 
of the operator $H$ is proved. A scattering theory is developed 
for $H$ in the case where the operator $A$ has continuous 
spectrum.
\end{Abstract}

\Section{Introduction}

Perturbations, depending on the spectral parameter (usually energy of
system) arise in a lot of quantum--mechanical problems
typically (see e.g. Ref.~[1]) as a result  of dividing the Hilbert
space $\cH$ of physical system in two subspaces,
$\cH=\cH_1\oplus\cH_2$.  The first one, say $\cH_1$, is interpreted
as a space of some ``external''  degrees of freedom. The second one,
$\cH_2$, is associated with an ``internal'' structure of the system.
The Hamiltonian ${\bH}$ of the system looks as a matrix,
\be
\label{twochannel}
{\bH}=\left[
  \begin{array}{cc}
                     A_1              &          B_{12}          \\
                     B_{21}           &          A_2
\end{array}
\right]
\ee
with $A\ad$, $\alpha=1,2$, the channel Hamiltonians (self-adjoint
operators in $\cH\ad$) and
$B_{12}$, $B_{21}=B_{12}^{*}$, the coupling operators.
Reducing the spectral problem
${\bH}U=zU$ , $U=\{ u_1 ,u_2 \}$ to the
channel $\an$ only one gets the spectral problem
\be
\label{ini}
[A\ad  +V\ad  (z)]u\ad  =zu\ad  ,
\quad \an =1,2,
\ee
where the perturbation
\be
\label{epot}
V\ad  (z)=-B\abd (A\bd -z I\bd)^{-1}B\bad ,
\quad \bn \neq \an,
\ee
depends on the spectral parameter $z$ as the resolvent
$(A\bd - z I\bd)^{-1}$ of the Hamiltonian $A\bd$. Here, by $I\bd$ we
understand the identity operator in $\cH\bd$.
\Newpara
The present paper is a summary of the author's
works~[2]---[4] considering a possibility to
``remove'' the energy dependence from  perturbations of the
type~(\ref{epot}). Namely, in~[2]---[4] we
search for such a new perturbation (``potential'') $W\ad$ not
depending on $z$ that spectrum of the Hamiltonian $H\ad=A\ad+W\ad$ is
(a part of) the spectrum of the problem~(\ref{ini}).  At the same
time, the respective eigenvectors of $H\ad$ become also those
for~(\ref{ini}).  An interest to the problem of such a removal of
dependence on the spectral parameter from perturbations is stimulated
in particular by a rather conceptual question (see for instance
Ref.~[5]) concerning a use of the
two--body energy--dependent potentials in few--body nonrelativistic
scattering problems. Since the energies of  pair subsystems are not
fixed in the N--body (N$\geq 3$) system, a direct embedding of such
potentials into the few--body Hamiltonian is impossible. Thus, the
replacements of the type (\ref{epot}) energy--dependent potentials
with the respective new potentials $W_\alpha$ could be considered as
a way to overcome this difficulty (see Ref.~[4] for
discussion).
\Newpara
The Hamiltonians $H\ad$ are found
in~[2]---[4] as solutions of the
non-linear operator equations (first appeared in Ref.~[6]) 
\be
\label{basic}
H\ad=A\ad + V\ad (H\ad).
\ee
The operator-value function $V\ad(Y)$ of the operator variable $Y$,
$Y:\, \cH\ad\rightarrow\cH\ad,$ is defined by us 
in such a way [see formula (\ref{Vdef})] that 
eigenfunctions $\psi$ of $Y$, $Y\psi =z\psi$,
become automatically those for $V\ad(Y)$ and
$V\ad(Y)\psi=V\ad(z)\psi.$ We have proved a solvability of
Eq.~(\ref{basic}) in the case where the Hilbert--Schmidt norm
$\|B\abd\|_2$ of the operators $B\abd$ satisfies the condition
$\|B\abd\|_2<\frac{1}{2}{\rm dist}\{\sigma(A_1),\sigma(A_2)\}$ in
supposition that spectra $\sigma(A\ad)$ of the operators $A\ad$ are
separated, ${\rm dist}\{\sigma(A_1),\sigma(A_2)\}>0$ (see Theorem~1).
\Newpara
In Ref.~[2], the problem of the removal of
energy dependence from the type (\ref{epot}) perturbations was
considered in details when one of the operators $A\ad$ is the
Schr\"{o}dinger operator in $L_2({\R}^n)$ and another one has a
discrete spectrum only. The report
[3]  announces the results concerning the
equations~(\ref{basic}) and properties of their solutions $H\ad$ in a
rather more general situation where the Hamiltonian ${\bH}$ may be
rewritten in terms of a two-channel variant of the Friedrichs model
investigated by {\sc O.A.Ladyzhenskaya} and {\sc L.D.Faddeev}~[7].
In particular in~[2] and~[3] a scattering
problem is studied for $H\ad$ in the case if $A\ad$ has continuous
spectrum and the basis property of the eigenfunction system of the
operator $H\ad$ is shown.
\Newpara
In the paper~[4], we specify the assertions
from~[3] and give proofs for them.
Also, we pay attention to an important circumstance disclosing a
nature of solutions of the basic equations (\ref{basic}). Thing is
that the operators $W\ad=V\ad(H\ad)$ may be present in the form
$W\ad=B\abd Q\bad$ with $Q\bad$ satisfying the stationary Riccati
equations~(\ref{QbasicSym}). Exactly the same equations always arise
if one makes a block diagonalization of
the type~(\ref{twochannel}) operator
matrices in the way described below in Lemma~1, so that the solutions
$H\ad$, $\alpha=1,2,$ of Eqs.~(\ref{basic}) determine in fact parts
of the operator ${\bf H}$ in respective invariant subspaces. The idea
of such a diagonalization was applied already by {\sc S.Okubo}~[8]
to some quantum--mechanical Hamiltonians. It was used later by
{\sc V.A.Malyshev} and {\sc R.A.Minlos}~[9] in a method of
construction of invariant subspaces for a class of self--adjoint
operators in statistical physics.  This idea was used also in the
recent paper~[10] by {\sc V.M.Adamjan} and {\sc H.Langer} who
studied spectral properties of a class of the type~(\ref{ini})
spectral problems and in particular, a possibility to choose among
their solutions a Riesz basis in $\cH\ad$.

\Section{Construction of the operators $H\ad$}

We study the spectral problem~(\ref{ini}) with  perturbation
$V\ad(z)$  given by~(\ref{epot}).  We suppose that $B\bad$ is a
linear operator from ${\cH}\ad$ to ${\cH}\bd$ with a finite
Hilbert--Schmidt norm $\|B\bad\|_2$, $\|B\bad\|_2<\infty$.  A goal of
the work is a construction  of  such an operator $H_{\alpha }$ that
its each eigenfunction $u_{\alpha }, H_{\alpha }u_{\alpha
}=zu_{\alpha }$, together with eigenvalue $z$,  satisfies
Eq.~(\ref{ini}).  The operator $H\ad$ is searched for as a solution
of the non-linear operator equation (\ref{basic}).  To obtain this
equation we introduce the following operator-value function $V\ad
(Y)$ of the  operator variable $Y:\quad$
\be
\label{Vdef}
V\ad  (Y)=B\abd\Int_{\sigma \bd }
E\bd (d\mu )B\bad (Y-\mu I\ad)^{-1},
\ee
$Y:$ ${\cH} \ad  \rightarrow {\cH} \ad  $ .
Here, $\sigma\bd$ is the spectrum and $E\bd$,
the spectral measure of the operator $A\bd$.
The integral over $E\bd$ in (\ref{Vdef})
for $Y$ such that $\Sup{\mu\in\sigma\bd}\|(Y-\mu I\ad)^{-1}\|<\infty$
may be constructed in the same way as the usual
integrals of scalar functions
over spectral measure.
For $\|B\bad\|_2<\infty$, the existence of this  integral
as a bounded operator from $\cH\ad$ to $\cH\bd$
is proved in~[4].
We notice that if $\phi$ is an eigenfunction of $Y$,
$Y\phi=z\phi$, then automatically
$
V\ad(Y)\phi=B\abd\Int_{\sigma\bd} E\bd(d\mu)
B\bad(z-\mu)^{-1}\phi =B\abd(z-A\bd )^{-1}B\bad\phi
=V\ad(z)\phi.
$
This means that
$H\ad$ satisfies with its eigenfunctions $\psi\ad$ the  relation
$H\ad\psi\ad=(A\ad+V\ad(H\ad))\psi\ad$
and one can spread this relation
over all the linear combinations of the eigenfunctions.
Supposing that the  eigenfunctions system  of $H\ad$  is
dense in $\cH\ad$ one spreads this
equation over all the domain ${\cal D}(A\ad)$.
As a result we come to the desired {\it basic equation} (\ref{basic})
for $H\ad$ (see also~[2]---[4] and Refs. therein).
Eq.~(\ref{basic}) means that the construction of the operator $H\ad$
is reduced to the searching for the operator
$Q\bad=\Int_{\sigma\bd} E\bd(d\mu)B\bad(H\ad-\mu I\ad)^{-1}.$
Since $H\ad=A\ad+B\abd Q\bad$, we have
\be
\label{Qbasic1}
Q\bad =
\Int_{\sigma \bd }E\bd (d\mu )
B\bad (A\ad  +
B\abd Q\bad -\mu I\ad)^{-1},
\quad \bn \neq \an.
\ee
We restrict ourselves to a study of
Eq.~(\ref{Qbasic1}) solvability  only  in  the  case where spectra
$\sigma_{1}$ and $\sigma_{2}$ are  separated,
$d_{0}={\rm dist}(\sigma _{1},\sigma _{2})>0.$
Applying to Eq.~(\ref{Qbasic1}) the contracting  mapping  theorem,
one comes to the following:

\theo{\Theorem}{1}
{Let $M\bad(\dn)$ be a set of bounded operators $X,$ $X:$
$\cH\ad\rightarrow\cH\bd$, satisfying the inequality
$\|{X}\|$ $\le\dn$ with $\dn>0$. If this $\dn $ and the norm
$\|{B\abd }\|_{2}$ satisfy the condition
$\|B\abd\|_{2}<d_{0}\Min{}\{\Frac{1}{1+\dn},\Frac{\dn }{1+\dn^2}\}$,
then Eq.~(\ref{Qbasic1}) is  uniquely solvable in $M\bad(\dn)$.
In particular the equation (\ref{Qbasic1}) is uniquely solvable
in the unit ball $M\bad(1)$ for any $B\abd$ such that
$\|B\abd\|_2<\Frac{1}{2}d_0.$}
Eq.~(\ref{Qbasic1}) can be rewritten
(see~[3], [4])
also in  symmetric form as a stationary Riccati equation,
\be
\label{QbasicSym}
Q\bad A\ad  -
A\bd Q\bad +
Q\bad B\abd Q\bad =B\bad .
\ee
One finds immediately from  Eqs.~(\ref{QbasicSym}), $\an=1,2$, that
if $Q\bad$ gives a solution $H\ad=A\ad+B\abd Q\bad$ of the problem
(\ref{basic}) in the channel $\an$  then
$Q\abd=-Q\bad^{*}=-\Int_{\sigma\ad }( H^{*}\ad -\mu I\ad)^{-1}B\abd
E\bd(d\mu )$ gives an analogous solution $H\bd=A\bd+B\bad Q\abd$ in the
channel $\bn$.

\theo{\Lemma}{1}
{Let $Q\bad$ and $Q\abd=-Q\bad^{*}$ be
solutions of Eqs.~(\ref{QbasicSym}).
Then the transform ${\bH}'={\cal Q}^{-1}{\bH}{\cal Q}$ with
$
{\cal Q}=\left[\begin{array}{lr}   I_1      &   Q_{12} \\
                                   Q_{21}   &     I_2
\end{array}\right]
$
reduces the operator ${\bH}$ to the block--diagonal form,
${\bH}'=\diag\{H_1,H_2\}$ with $H\ad=A\ad+B\abd Q\bad.$}
One can find assertions analogous to Lemma~{1}
in Refs.~[9] and~[10].
A solvability (for sufficiently small $\|B\abd\|$)
of the equation~(\ref{QbasicSym}) was proved
in~[9], [10] for different
situations and by rather different methods but
also in the supposition ${\rm dist}\{\sigma(A_1),\sigma(A_2)\}>0$.

\theo{\Remark}{1}{\rm Let $X\ad=I\ad-Q\abd Q\bad=I\ad+Q\abd Q\abd^{*}$.
It follows from Lemma~{1} that
the operator $\tilde{\cQ}=\cQ X^{-1/2} $ with $X=\diag\{X_1,X_2\}$ is
unitary.  Thus, the operator
${\bH}''=\tilde{\cQ}^{*}\bH\tilde{\cQ}=X^{1/2}{\bH}'X^{-1/2}$ becomes
self-adjoint in $\cH $.  Since $ {\bH}''=\diag\{ H''_1, H''_1 \}$
with $ H''\ad=X\ad^{1/2} H\ad X\ad^{-1/2}$, the operators $H''\ad$,
 $\,\,\,\an=1,2,$ are self-adjoint on $\cD(A\ad)$ in $\cH\ad $.
Moreover, the operators
${\bH}\au=\tilde{\cQ}\cdot\diag\{H''\ad,0\}\cdot\tilde{\cQ}^{*}=$
$\cQ\cdot\diag\{H\ad,0\}\cdot\cQ^{-1}$ represent parts of the
Hamiltonian ${\bH}$ in the corresponding invariant subspaces
$\cH\au=\{f:\,\,f=\{f\ad,f\bd\}\in\cH,\,\, f\ad\in\cH\ad,\,\,
f\bd=Q\bad f\ad \}$ (see also Refs.~[9],
[10]).}

\Section{Spectra of the Hamiltonians $H\ad$ and basis
properties of their eigenfunctions}

Let us suppose that $Q\bad$ and $Q\abd=-Q\bad^{*}$ are solutions of
Eqs.~(\ref{Qbasic1}) and~(\ref{QbasicSym}) which are spoken about in
Theorem~1. Since we take $B\abd$ with
$\|B\abd\|\leq \|B\abd\|_2<d_0/2$ and $\|Q\bad\|<1$,
the spectra $\sigma(H_1)$ and $\sigma(H_2)$ do not intersect
(actually, when these spectra are discussed in
Refs.~[3], [4], a more general case is also
considered where not necessary $\|Q\bad\|<1$).
 By Lemma~1, the operator
${\bH}'=\diag\{ H_1, H_2\}$ is connected with
the (self-adjoint) operator  ${\bH}$
by a similarity transform. Thus, the spectra
$\sigma(H_1)$ and $\sigma(H_2)$ of the operators $H\ad$, $\an=1,2$,
are real and $\sigma(H_1)\bigcup\sigma(H_2)=\sigma({\bH})$.
Continuous spectrum $\sigma_c (H\ad)$ of every $H\ad$
coincides with that, $\sigma\ad^c$, of the operator $A\ad$,
$\sigma_c (H\ad)=\sigma\ad^c$, since due to $\|B\abd\|_2 < +\infty$,
the potential $W\ad=B\abd Q\bad $ is a compact operator.
\Newpara
For more concrete statements concerning the spectra of the operators
$H\ad$ we accept some presuppositions
restricting us as regards ${\bf H}$
to the case of a two--channel variant of the Friedrichs model in the
form~[7] reproducing often encountered
quantum--mechanical situations. At first, we assume that the operator
${\bf H}$ is defined in that representation where the operators
$A\ad,$ $\an =1,2,$ are diagonal. We suppose
that the continuous  spectra $\sigma^{c}\ad$
are absolutely continuous and consist of a finite
number of finite (and may be one or two infinite) intervals.
At second, we suppose that discrete spectra $\sigma\ad^d$ of
the operators $A\ad$, $\an=1,2,$ do not intersect with $\sigma\ad^c$,
$\sigma\ad^d \bigcap\sigma\ad^c =\emptyset$,
and consist of a finite number of points with finite multiplicity.
The coupling operators $B\abd$ are supposed to be the integral ones
with sufficiently quickly decreasing (in the case of unbounded
$\sigma\ad^c$) kernels being smooth in
the H\"older sense (see Refs.~[3], [4]
for details).
\Newpara
With these presuppositions the continuous spectrum
$\sigma_{c}({\bH})=\sigma ^{c}_{1}\cup \sigma ^{c}_{2}$  of the
operator ${\bf H}$ is absolutely continuous and its part ${\bH}^{c}$
acting in respective invariant  subspace,
is unitary equivalent to the operator
${\bH}_{0}=A^{(0)}_{1}\oplus A^{(0)}_{2}$ with $A^{(0)}\ad,$
$\an=1,2,$ the part of $A\ad$ acting in the invariant
subspace $\cH\ad^c$ corresponding to $\sigma\ad^c$.
Namely, there exist the wave operators $U^{(+)}$ and $U^{(-)},$
$U^{(\pm )}=\left(\begin{array}{lr}
u^{(\pm )}_{11}    &     u^{(\pm )}_{12}  \\
u^{(\pm )}_{21}    &     u^{(\pm )}_{22}
\end{array}\right)=$ $s\!-\!\Lim{t\rightarrow\mp\infty}
{\rm e}^{i{\bH}t}{\rm e}^{-i{\bH}_0 t},$
with the properties:
${\bH}U^{(\pm )}=U^{(\pm )}{\bH}_{0},$ $U^{(\pm )*}U^{(\pm )}=I,$
$U^{(\pm )}U^{(\pm )*}=I-P_d.$
Here, by $P_d$ we understand the orthogonal projector on the
subspace corresponding to the discrete spectrum $\sigma _{d}({\bH})$
of ${\bH}$. The kernel $u^{(\pm )}_{\an\an}(\lambda,\lambda')$ of the
component $u^{(\pm)}_{\an\an},\an=1,2,$ represents
a (generalized) eigenfunction of continuous
spectrum of the problem~(\ref{ini})
for $z=\lambda'\pm i0,$ $\lambda'\in \sigma^{c}\ad$.  At the same
time $u^{(\pm )}\abd (\lambda,\lambda')$ is the problem~(\ref{ini})
eigenfunction corresponding  to $\lambda '\in \sigma\bd^c$.
\Newpara
By $U_{j},$ $j=1,2,\ldots$, we denote eigenvectors,
$U_{j}=\{u\ju_{1},u\ju_{2}\},$ $\|U_{j}\|=1,$ and by $z_{j},$ $
z_{j}\in\R$, the respective eigenvalues of
$\sigma_d({\bH})$. We assume that in the case of multiple
discrete eigenvalues, certain $z_j$ may be repeated in the
numeration. The component $u\ju\ad$ of the vector
$U_{j}$ is a solution of Eq.~(\ref{ini}) for $z=z_j$.
\Newpara
Let us return, with the presuppositions above, to the operators $H\ad$.
First, let us assume $\sigma _{d}(H\ad)\neq\emptyset$. Then, it
follows from the construction  of the function (\ref{Vdef}) that if
$z\in\sigma_{d}(H\ad)$  then this $z$  becomes automatically a point of
the discrete spectrum of the initial spectral problem~(\ref{ini}).
At the same time $\psi\ad $ becomes its eigenfunction.  We shall
denote the eigenfunctions of the $H\ad$ by $\psi\ju\ad,$
$\psi\ju\ad=u\ju\ad$, keeping for them the same numeration as for
the eigenvectors $U_j$, $U_j=\{u\ad\ju,u\bd\ju\},$ of the
Hamiltonian ${\bH}$, ${\bH}U_j =z_j U_j$, $z_j\in\sigma_d ({\bH}).$
Respective eigenvectors of the adjoint operator $H^*\ad$,
$H^{*}\ad=A\ad+Q^{*}\bad B\bad,$ are
$\tilde{\psi}\ju\ad=\psi\ju\ad-Q\abd{u}\ju\bd$.
Due to Lemma~{1},
$\sigma_d({\bf H})=\sigma_d(H_1)\bigcup\sigma_d(H_2).$
Since in conditions of Theorem~1
$\sigma(H_1)\bigcap\sigma(H_2)=\emptyset$, we have also
$\sigma_d(H_1)\bigcap\sigma_d(H_2)=\emptyset$.
\Newpara
Dealing with the continuous spectrum of $H\ad$ we take into account
the fact that the solutions $Q\bad$ and $Q\abd=-Q\bad^{*}$ of
Eqs.~(\ref{Qbasic1}) and~(\ref{QbasicSym}) corresponding to the
operators $B\abd$ with H\"older kernels, have the H\"older kernels
themselves. The same is true as well for $W\ad=B\abd Q\bad$. Then we
can prove~[3], [4] that the operators
$\Psi\ad^{(\pm)}=u_{\an\an}^{(\pm)}$ turn out to be the wave
operators between $H\ad$ and $A\ad^{(0)}$:
$\Psi\ad^{(\pm)}=s\!-\!\Lim{t\rightarrow\mp\infty}
\exp(iH\ad t)\exp(-iA\ad^{(0)}t)$ and
$H\ad\Psi^{(\pm)}=\Psi\ad^{(\pm)}A\ad^{(0)}$.  At the same time the
operators
$\tilde{\Psi}^{(\pm)}\ad=\Psi^{(\pm)}\ad-{Q\abd}u^{(\pm)}\bad$
become those ones for $H\ad^*$.

\theo{\Theorem}{2}{
The following orthogonality relations take place:
$\langle\psi\ju\ad,\tilde{\psi}^{(k)}\ad\rangle=\dn_{jk},$
$\Psi^{(\pm)*}\ad\tilde{\Psi}^{(\pm)}\ad=\reduction{I\ad}{\cH\ad^c},$
$\tilde{\Psi}^{(\pm )*}\ad\psi\ju\ad=0$ and
$\Psi\ad^{(\pm)*}\tilde{\psi}\ju\ad=0.$ Also,
the completeness relations are valid,
$\Sum_{j: H\ad u\ju\ad=z_j u\ju\ad}
\psi\ju\ad\langle\cdot,\tilde{\psi}\ju\ad\rangle+
\Psi^{(\pm)}\ad\tilde{\Psi}^{(\pm)*}\ad=I\ad,\quad\an=1,2.$
For all this $ S^{(\an)}=$ $ \Psi\ad^{(-)-1}\Psi\ad^{(+)}=$ $
\tilde{\Psi}\ad^{(-)*}\Psi\ad^{(+)}=$ $ \Psi\ad^{(-)*} X\ad
\Psi\ad^{(+)}$ represents a scattering operator for a system
described by the Hamiltonian $H\ad$. In fact, this operator coincides
with the component $S\aad$ of the scattering operator $S$,
$S=U^{(-)*}U^{(+)}$, for a system described by the two--channel
Hamiltonian ${\bH}$.}

\begin{Acknowledgements}
The author is thankful to Prof.~R.Mennicken for his interest
to this work and to Professors  V.M.Adamjan, H.Langer, R.A.Minlos 
and A.A.Shkalikov for valuable discussions. Especially, the author 
is grateful to Prof.~B.S.Pavlov and Dr.~K.A.Makarov for 
useful remarks and stimulating support.
\end{Acknowledgements}

\begin{References}
\item {\normalsize\sc  R.L.Jaffe, R.L.,  Low, F.E.}: Connection between
        quark--model eigenstates and low--energy scattering;  Phys. Rev.
         {\bf D19} (1979), 2105--2118; {\normalsize\sc Simonov, Yu.A.}:
           Hadron--hadron interaction in the compound--bag model;
           Yadernaya Fiz. {\bf 36} (1982), 722--731 [Russian].%
\item   {\normalsize\sc Motovilov, A.K.}:
        The removal of an energy dependence from
        the interaction in two--body systems; J.Math.Phys. {\bf 32} (1991),
        3509--3518.%
\item  {\normalsize\sc Motovilov, A.K.}: Potentials appearing after the
        removal of an energy dependence and scattering by them;
        Proc. Intern. Workshop ``Mathematical
        Aspects of the Scattering Theory and Applications''
        (St.Petersburg, 1991),~101--108.%
\item   {\normalsize\sc Motovilov, A.K.}:
        Removal of the energy dependence from
        the resolvent--like energy--dependent interactions;
        Preprint JINR E5--94--259 [nucl--th/9505030] 
        (available via e--mail from
        Los--Alamos e--print database by sending the empty message to
        e--address {\tt nucl-th@xxx.lanl.gov} with
        Subject: {\tt get 9505030}); to appear in
        Teor. Mat. Fiz. {\bf 104} (1995).%
\item
       {\normalsize\sc McKellar, B.H.J., McKay C.M.}:
        Formal scattering theory for energy--dependent potentials;
       { Aust.J.Phys.}{\bf 36} (1983).,
        607--616;
       {\normalsize\sc Schmid, E.W.}:
       The problem of using energy--dependent nucleon--nucleon 
       potentials in nuclear physics;  
       Helv. Phys.Acta {\bf 60} (1987), 394--397.%
\item   {\normalsize\sc Braun, M.A.}:
        On relation between quasipotential equation 
        and Schr\"odinger equation; Teor. Mat. Fiz. {\bf 72} (1987),
        394--402 [Russian].%
\item  {\normalsize\sc Ladyzhenskaya, O.A., Faddeev, L.D.}:
        On perturbation theory of continuous spectrum;        
        Doklady AN SSSR
        {\bf 120} (1958),  1187--1190 [Russian]; 
        {\normalsize\sc Faddeev, L.D.}:
	On the Friedrichs model in 
        perturbation theory of continuous spectrum;        
        Trudy Mat. In-ta AN SSSR {\bf 73} (1964), 292--313 [Russian].
\item {\normalsize\sc Okubo, S.}: Diagonalization of Hamiltonian
         and Tamm--Dancoff equation; Progr. Theor. Phys. {\bf 12} (1954),
         603--622.
\item  {\normalsize\sc Malyshev, V.A., Minlos, R.A.}:
       Invariant subspaces of clustering operators. I.,
       J. Stat. Phys. {\bf 21} (1979), 231--242;
       Invariant subspaces of clustering operators. II.,
       Comm. Math. Phys. {\bf 82} (1981), 211--226.%
\item {\normalsize\sc Adamjan, V.M., Langer, H.}:
          {Spectral properties of a class
           of rational operator-value functions}; to appear
           in {J.Operator Theory}.%
\end{References}
\begin{Address}
{\sc Dr. A.K.Motovilov,}
Laboratory of Theoretical Physics, JINR,
141980 Dubna, Moscow region, Russia\\
E--address: motovilv@thsun1.jinr.dubna.su
\end{Address}

\end{document}